\documentclass[11pt]{amsart}
\usepackage{amsmath}
\usepackage{amssymb}
\usepackage{amsfonts}
\usepackage{latexsym}
\usepackage{color}

\newtheorem{defn}{Definition}[section]
\newtheorem{thm}[defn]{Theorem}
\newtheorem{prop}[defn]{Proposition}
 \newtheorem{lemma}[defn]{Lemma}
\newtheorem{cor}[defn]{Corollary}
 \newtheorem{example}[defn]{Example}
\newtheorem{rem}[defn]{Remark}
\newtheorem{rems}[defn]{Remarks}
\newcommand{\belem}{\begin{lemma}}
\newcommand{\enlem}{\end{lemma}}

\newcommand{\beprop}{\begin{prop}}
\newcommand{\enprop}{\end{prop}}

\newcommand{\betheo}{\begin{theorem}}
\newcommand{\entheo}{\end{theorem}}

\newcommand{\becor}{\begin{cor}}
\newcommand{\encor}{\end{cor}}

\def\x{\relax\ifmmode {\mbox{*}}\else*\fi}
\newcommand{\beex}{\begin{example}$\!\!${\bf }$\;$\rm }
	\newcommand{\enex}{ \end{example}}

\newcommand{\berem}{\begin{rem}$\!\!${\bf }$\;$\rm }
	\newcommand{\enrem}{ \end{rem}}
 \newcommand{\berems}{\begin{rems}$\!\!${\bf }$\;$\rm }
 	\newcommand{\enrems}{ \end{rems}}

\newcommand{\bedefi}{\begin{defn}$\!\!${\bf }$\;$\rm }
	\newcommand{\findefi}{\end{defn}}

\newcommand{\be}{\begin{equation}}
\newcommand{\en}{\end{equation}}
\newcommand{\bea}{\begin{eqnarray}}
\newcommand{\ena}{\end{eqnarray}}
\newcommand{\beano}{\begin{eqnarray*}}
	\newcommand{\enano}{\end{eqnarray*}}
\newcommand{\bee}{\begin{enumerate}}
	\newcommand{\ene}{\end{enumerate}}
\newcommand{\bei}{\begin{itemize}}
	\newcommand{\eni}{\end{itemize}}

\newcommand{\norm}[2]{
\left\| #2 \right\|_{#1}}
\newcommand{\hs}{Hilbert space}
 \newcommand{\ov}{\overline}
\newcommand{\pip}{PIP-space}

\newcommand{\dis}{\displaystyle}

\newcommand{\mc}{\mathcal}

\newcommand{\NN}{\mathbb{N}}
\newcommand{\RN}{\mathbb{R}}
 \newcommand{\ZN}{\mathbb{Z}}
\newcommand{\CN}{{\mathbb C}}
  
 \newcommand{\bdim}{{\bf Proof. }}
 \newcommand{\edim}{\hfill $\square$\bigskip}

\newcommand{\Hil}{{\mc H}}

 \newcommand{\B}{{\mc B}}

 \newcommand{\D}{{\mc D}}
 \newcommand{\G}{{\mc G}}
 \newcommand{\cP}{{\mc P}}
 \newcommand{\T}{{\mc T}}

\newcommand{\ha}{^{\ast}}
\newcommand{\ip}[2]{\langle {#1}|{#2}\rangle}

\def\H{{\mathcal H}}

\newcommand{\ud}{\,\mathrm{d}}
\newcommand{\un}{\underline}
	\definecolor{ao}{rgb}{0.0, 0.5, 0.0}
\numberwithin{equation}{section}

\begin{document}

 \title{Weak $A$-frames and weak $A$-semi-frames}

\author {Jean-Pierre Antoine   $^1$ }       
\author {Giorgia Bellomonte  $^2$   }
\author {Camillo Trapani $^3$} 
\address{Jean-Pierre Antoine   \newline
Universit\'e catholique de Louvain \newline
Institut de Recherche en Math\'ematique et  Physique\newline
 B-1348   Louvain-la-Neuve, Belgium}
 \email {jean-pierre.antoine@uclouvain.be}

\address{Giorgia Bellomonte    \newline 
 Universit\`a di Palermo\newline
Dipartimento di Matematica e Informatica \newline
I-90123 Palermo, Italy}
 \email {giorgia.bellomonte@unipa.it}
 
\address{Camillo Trapani  newline
Universit\`a di Palermo\newline
Dipartimento di Matematica e Informatica \newline
I-90123 Palermo, Italy}
 \email {camillo.trapani@unipa.it}
 
	\subjclass[2010]{41A99, 42C15}
\keywords{\em $A$-frames, {weak} (upper and lower) $A$-semi-frames,  lower atomic systems, $G$-duality}


 			\vspace{5pt}
		\noindent
   
\maketitle

 { \bf  Abstract}
After reviewing the interplay between frames and lower semi-frames,
we introduce the notion of  lower semi-frame controlled by a densely defined  operator $A$ or, for short, a
\emph{weak lower $A$-semi-frame} and we study its properties. In particular, we compare it with that of  lower atomic systems, introduced in \cite{GB}.
We discuss  duality properties and we suggest  several possible  definitions for weak $A$-upper semi-frames.
Concrete examples are presented. 

\textbf{Keywords: }{$A$-frames, {weak} (upper and lower) $A$-semi-frames,  lower atomic systems, $G$-duality}
		
		\vspace{5pt}
		\textbf{2010 Mathematics Subject Classification:} 41A99, 42C15.

 {\bf Corresponding author:} {Giorgia Bellomonte}
 
 
\medskip


\section{Introduction and basic facts}\label{sec-intro}

We consider an infinite dimensional Hilbert space $\H$  with inner product $\ip{\cdot}{\cdot}$, linear in the first entry, and norm $\|\cdot\|$. $GL(\H)$ denotes the set of all invertible bounded operators on $\H$ with bounded inverse.
Given a linear operator $A$, we denote its domain by $\mathcal{D}(A)$, its
range by $\mathcal{R}(A)$  and its adjoint by $A^*$, if $A$ is densely defined.
Given a locally compact, $\sigma$-compact
space $(X,\mu)$ with a   (Radon) measure $\mu$, a function $\psi:X\mapsto \H, x\mapsto \psi_x$ is said to be {\it weakly measurable}  if for every $f\in \H$ the function $x\mapsto \ip{f}{\psi_x}$ is measurable. As a particular case, we obtain a discrete situation if $X=\NN$ and $\mu$ is the counting measure.

Given a weakly measurable function  $\psi$, the operator $C_\psi:\D(C_\psi)\subseteq \H \to L^2(X,\ud\mu)$ with domain
$$
\D(C_\psi):=\left \{f\in \H: \int_{X}  |\ip{f}{\psi_{x}}| ^2 \, \ud \mu(x) <\infty \right\}
$$
and $({C}_{\psi}f)(x) =\ip{f} {\psi_{x}}, f \in \D(C_\psi),C_\psi$
is called the {\it analysis} operator of $\psi$.
	\berem
	In general the domain of $C_\psi $ is not dense, hence $C_\psi ^*$ is not well-defined. An example of function whose analysis operator is densely defined can be found in  \cite[Example 2.8]{GB}, where $\D(C_\psi )$ coincides with the domain of a densely defined sesquilinear form associated to $\psi$. Moreover, a sufficient condition for $\D(C_\psi )$ to be dense in $\H$ is that $\psi_x\in\D(C_\psi ) $ for every $x\in X$, see  \cite[Lemma 2.3]{AB}.
\enrem
	\begin{prop}\cite[Lemma 2.1]{AB}\label{prop: closedness C}
	Let $(X,\mu)$ be a  locally compact, $\sigma$-compact
		space, with a Radon measure $\mu$  and $\psi:x\in X\mapsto\psi_x\in\H$ a weakly measurable function. Then the analysis  operator $C_\psi $  is  closed.\end{prop}

Consider the set $ \D(\Omega_\psi) =\D(C_\psi)$
and the mapping
$\Omega_\psi : \D(C_\psi) \times\D(C_\psi) \to \CN$ defined by
\begin{equation}\label{defn}
\Omega_\psi (f, g) :=\int_X\ip{f}{\psi_x}\ip{\psi_x}{g}\ud\mu(x).
\end{equation}
$\Omega_\psi $ is clearly a nonnegative symmetric sesquilinear form which is  well defined for every $f,g\in\D(C_\psi)$ because of the Cauchy-Schwarz inequality. It is unbounded in general. Moreover, since $\D(C_\psi)$ is the largest domain such that $\Omega_\psi$ is defined on $\D(C_\psi)\times\D(C_\psi)$, it follows that
\begin{equation}\label{eq: Omega closed}
\Omega_\psi (f,g)=\ip{C_\psi f}{C_\psi g}, \quad\forall f,g\in\D(C_\psi ),
\end{equation}
where $C_\psi $ is the analysis operator defined above.  Since $C_\psi $ is a closed operator,
the form $\Omega_\psi$ is closed, see e.g. \cite[Example VI.1.13]{Kato}.
If $\D(C_\psi)$ is dense in $\H$, then by Kato's first representation
 theorem \cite[Theorem VI.2.1]{Kato} there exists a positive  self-adjoint  operator ${\sf T}_\psi$ {\em associated to the sesquilinear form} $\Omega_\psi$  on
\be\label{eq: Omega}
\D({\sf T}_\psi)  =\left\{f\in \D(\Omega_\psi): h \mapsto \int_X\ip{f}{\psi_x} \ip{\psi_x}{h}\ud\mu(x)
\mbox{ is bounded in } \D(C_\psi)\right\}
\en
defined  by
\be\label{def: Frame opera}
{\sf T}_\psi f:=h  
\en
with $h$ as in \eqref{eq: Omega}. The density of $\D(\Omega_\psi)$ ensures the uniqueness of the vector $h$.
The operator ${\sf T}_\psi$ is the  greatest one whose domain is contained in $\D(\Omega_\psi)$ and  such that
$$
\Omega_\psi(f, g) = \ip{{\sf T}_\psi f}{g},\qquad  f \in \D({\sf T}_\psi),\, g\in\D(\Omega_\psi).
$$
The set   $\D({\sf T}_\psi)$ is dense in\ $\D(\Omega_\psi)$  , see  \cite[p. 279]{Kato}. { In addition}, by Kato's second representation theorem \cite[Theorem VI.2.23]{Kato},  { we have } $\D(\Omega_\psi)=\D({\sf T}_\psi^{1/2})$
and $$
\Omega_\psi(f, g) = \ip{{\sf T}_\psi^{1/2}f}{{\sf T}_\psi^{1/2}g},\qquad \forall  f, g \in \D(\Omega_\psi),
$$
 { hence}, comparing with \eqref{eq: Omega closed}, we  { deduce} ${\sf T}_\psi=C_\psi ^*C_\psi =|C_\psi |^2$ on
$\D({\sf T}_\psi)$.

\bedefi
The operator ${\sf T}_\psi:\D({\sf T}_\psi) \subset\H\to\H$ defined by \eqref{def: Frame opera} will be called the {\em generalized frame operator} of the function $\psi:x\in X\to\psi_x\in\H$. \findefi

{ Now we recall a series of notions well-known in literature, see e.g. \cite{AAG,AB,  Kaiser}.}
A weakly measurable function $\psi$ is said to be  \begin{itemize}
	\item $\mu$-total if $\ip{f}{\psi_x} = 0$ for
	a.e. $x\in  X$ implies that $f = 0$;
\item a  \textit{continuous frame} of $\H$ if
	there exist constants $0 < {\sf m} \leq {\sf M} < \infty$ (the frame bounds) such that
$$
	{\sf m}\|f\|^2 \leq
	\int_{X}  |\ip{f}{\psi_{x}}| ^2 \, \ud \mu(x)   \leq { \sf M}  \norm{}{f}^2 , \qquad \forall \, f \in \H; \,$$

	\item  a \textit{Bessel mapping}  of $\H$ if there exists ${ \sf M} >0$ such that
$$\int_{X}  |\ip{f}{\psi_{x}}| ^2 \, \ud \mu(x)   \leq { \sf M}  \norm{}{f}^2 , \qquad \forall \, f \in \H; \,
$$
\item  an \emph{upper semi-frame} of $\H$ if
there exists ${ \sf M} <\infty$ such that
$$
0<\int_{X}  |\ip{f}{\psi_{x}}| ^2 \, \ud \mu(x)   \leq { \sf M}  \norm{}{f}^2 , \qquad \forall \, f \in \H, \, f\neq0
$$
i.e.  if  
it is  a $\mu$-total Bessel mapping;
\item a \emph{lower semi-frame}  of $\H$ if  there exists   a constant ${\sf m}>0$ such that
\begin{equation}\label{eq:lowersf}
{\sf m}  \norm{}{f}^2 \leq  \int_{X}  |\ip{f}{\psi_{x}}| ^2 \, \ud \mu(x) , \qquad \forall \, f \in \H.
\end{equation}\end{itemize}
Note that the integral on the right hand side in \eqref{eq:lowersf} may diverge for some $f\in\H$, namely, for $f\not\in \D(C_\psi)$. Moreover, if $\psi$ satisfies \eqref{eq:lowersf} then it is automatically $\mu$-total.

\section{From semi-frames to frames and back}\label{sec2}

Starting from a  lower semi-frame, one can easily obtain a genuine frame, albeit in a smaller space. 
 Indeed, we have proved a theorem \cite[Prop.3.5]{ACT}, which implies  the following :

\beprop \label{prop5}
A weakly measurable function $\phi$ on $\H$  is a lower semi-frame of $\H$ whenever
 $\D(C_\phi)$ is complete for the norm $\norm{C_\phi}{f}^2=
 \int_{X}  |\ip{f}{\phi_{x}}| ^2 \, \ud \mu(x) = \|C_\phi f\|^2$,
continuously embedded into $\H$ and for some  $\alpha,{\sf m}, {\sf M} > 0$, one has
\bea
&&\alpha \norm{}{f} \leq \norm{C_\phi}{f} \; \text{and}  \label{eq:triple1}
\\
&&{\sf m}  \norm{C_\phi}{f}^2  \leq   \int_{X}  |\ip{f}{\phi_{x}}| ^2 \, \ud \mu(x)  \leq {\sf M}  \norm{C_\phi}{f}^2 ,  \; \forall \, f \in \D(C_\phi). \label{eq:triple2}
\ena
\enprop
Note that \eqref{eq:triple2} is trivial here.

 Following
the notation of  our previous papers,  denote by $\H( {\sf T}_\phi^{1/2})$ the \hs\  $ \D({\sf T}_\phi^{1/2})$  with the   norm  $\norm{1/2}{f}^2=
\norm{}{{\sf T}_\phi^{1/2}  f}^2$, where ${\sf T}_\phi$ is the generalized frame operator defined in \eqref{def: Frame opera}.
 In the same way, denote by  $\H( C_\phi)$ the \hs\  $ \D(C_\phi)$  with the
 inner product  $\ip{\cdot}{\cdot}_{C_\phi} =\ip{C_\phi\cdot}{C_\phi\cdot}$,
 and the corresponding norm $\norm{C_\phi  }{f}^2= \norm{}{C_\phi f}^2$. Then clearly   $\H( {\sf T}_\phi^{1/2}) = \H(C_\phi)$.

What we have obtained in Proposition \ref{prop5} is a frame in $\H(C_\phi)=\H( {\sf T}_\phi^{1/2}) $.
Indeed assume that  $\D(C_\phi)$ is dense.
Then, for every $x\in X$, the map $f\mapsto \ip{f}{\phi_x}$ is a bounded linear functional on the \hs\    $ \H(C_\phi)$.
 By the Riesz Lemma, there exists  an element $\chi_x^{\phi} \in \D(C_\phi)$ such that
\begin{equation*} 
\ip{f}{\phi_x} = \ip{f}{\chi_x^{\phi}}_{C_\phi}  \quad \forall  \, f \in \D(C_\phi).
\end{equation*}
  By Proposition \ref{prop5}, $\chi^{\phi}$ is a frame.

{Actually, one can say more \cite{ACT}. The   norm   $\norm{1/2}{f}^2= \norm{}{{\sf T}_\phi^{1/2}  f}^2$, is equivalent to the the graph norm  of ${\sf T}_\phi^{1/2}$. Hence $\ip{f}{\phi_x} =\ip{f}{\chi_x^{\phi}}_{C_\phi}=\ip{ f}{{\sf T}_\phi\chi_x^{\phi}}$ for all $f\in \D(C_\phi)$. Thus $\chi_x^{\phi}={\sf T}_\phi^{-1} \phi_x$ for all $x\in X$, i.e. $\chi^\phi$  is the canonical dual Bessel mapping of $\phi$ (we recall that  $\phi$ may have several duals).}

 \beprop
 	Let $\phi$ be a lower semi-frame of $\H$ with $\D(C_\phi)$ dense. Then the  canonical dual Bessel mapping  of $\phi$ is a tight frame for the Hilbert space 	$\H(C_\phi)$.
 \enprop

Conversely,    starting with a frame $\chi \in \D(C_\phi)$, does there exist a lower semi-frame $\eta$  of $\H$ such that $\chi $ is the frame $\chi^{\eta}$ constructed from $\eta$ in the  way described above. The answer is formulated in the following \cite[Prop. 6]{corso}
\beprop
Let $\chi $ be a frame of $\H(C_\phi)=\H({\sf T}_\phi^{1/2} )$. Then
\\
(i) there exists a lower semi-frame $\eta$ of $\H $ such that $\chi  = \chi^{\eta}$  if, and only if, $\chi  \in \D({\sf T}_\phi)$;
\\
(ii) if $\chi  =\chi^{\eta}$ for some lower semi-frame $\eta$ of $\H $, then $\eta= {\sf T}_\phi\chi$.
\enprop
\medskip

So far we have discussed the interplay between frames and lower semi-frames. But one question remains: how does one obtain semi-frames?  A standard construction is to start from an  unbounded operator $A$ and build a lattice of \hs s out of it, as described in \cite{AT-PIPmetric} and in \cite{ACT}. As we will see in Section \ref{sec-examples} (1) and (2) below, this approach indeed generates  a weak lower $A$-semi-frame.

 Before that, we need a new ingredient, namely the notion of metric operator.

Given    a closed unbounded    operator  $S$ with dense domain  $\D(S)$, define  the  operator $G= I+ S^\ast S$, which is unbounded, with $G > 1$ and bounded inverse. This is a \emph{metric operator}, that is,   a strictly positive
self-adjoint operator $G$, that is, $G > 0$ or $\ip{Gf}{f} \geq 0$ for every $f \in \D(G)$ and $\ip{Gf}{f} =
0$  if and only if $f= 0$.

Then the norm
$\|f\|_{G^{1/2}}=\|G^{1/2}f\|$
 is equivalent to  the graph norm of  $G^{1/2}$ on $\D(G^{1/2}) = \D(S)$ and
makes the latter into a \hs\ continuously embedded into  $\H$, denoted   by $\H(G) $.
Then $\H(G^{-1})$, built in the same way from $G^{-1}$,
 coincides, as a vector space, with the conjugate dual of  $\H(G)$.
 On the other hand,  $G^{-1}$ is bounded.
Hence we get the triplet
\be\label{eq:tri>1}
\H(G) \; \subset\; \H \; \subset\; \H(G^{-1}) = \H({G})^\times.
\en

Two developments arise from these relations. First, the triplet \eqref{eq:tri>1}
 is the central part of the discrete scale of Hilbert spaces $V_{\G}$  built on the powers of  ${G^{1/2}}$.
This means that $V_{\G}:= \{\H_{n}, n \in \ZN \}$,
where $\H_{n} = \D(G ^{n/2}),  n\in \NN$, with a norm equivalent to the graph norm, and $ \H_{-n} =\H_{n}^\times$:
\begin{equation*} 
 \ldots\subset\; \H_{2}\; \subset\;\H_{1}\; \subset\; \H \; \subset\; \H_{-1} \subset\; \H_{-2} \subset\; \ldots
\end{equation*}
Thus $\H_{1} =  \H({G}^{1/2} ) = \D(S)$,  $\H_{2} =  \H({G} ) = \D(S^\ast S)$, and $\H_{-2} =  \H({G}^{-1})$,
and so on. What we have obtained in this way is a Lattice of Hilbert Spaces (LHS), the simplest example of a
Partial Inner Product Spaces (\pip). See our monograph \cite{pip-book} about this structure.

One  may also add the end spaces of the scale, namely,
\be \label{eq:endscale}
\H_{\infty}({G} ):=\cap_{n\in \ZN} \H_n, \qquad \H_{-\infty}({G} ):=\bigcup_{n\in \ZN} \H_{n}.
\end{equation}
In this way, we get a genuine Rigged Hilbert Space:
\begin{equation*}
\H_{\infty}(G ) \subset \H \subset \H_{-\infty}({G} ).
\end{equation*}
In fact, one can go one more step  farther. Namely, following \cite[Sec. 5.1.2]{pip-book}, we can use quadratic interpolation theory \cite{berghlof} and build a continuous scale of Hilbert spaces
$\H_{\alpha},\alpha\geq 0$,
where $\H_{\alpha}=  \D(G^{\alpha/2})$,  {with the graph norm  $\|\xi\|_{\alpha}^2 = \|\xi\|^2 + \|G^{\alpha/2}\xi\|^2$ or, equivalently, the norm
$\norm{}{(I+G)^{\alpha/2}\xi}^2$.
Indeed every $G^\alpha, \alpha\geq 0$, is   an unbounded metric operator.
Next we define $\H_{-\alpha} =\H_{\alpha}^\times$ and thus obtain the full continuous scale $V_{\widetilde \G}:= \{\H_{\alpha}, \alpha \in \RN \}$.
Of course, one can replace $\ZN$ by $\RN$ in the definition \eqref{eq:endscale} of the end spaces of the scale.
\medskip

{A second development of the previous analysis
is that we have made a link to  the formalism based on metric operators that we have developed for the theory of pseudo-Hermitian operators, in particular non-self-adjoint Hamiltonians, as encountered in the so-called pseudo-Hermitian or $\cP\T$-symmetric quantum mechanics.
This is not the place, however, to go into details, instead we refer the reader to \cite{AT-PIPmetric,AT-metric2}  for a complete mathematical treatment.}

\section{Weak lower $A$-semi-frames}\label{sec3}

The following concept was introduced and studied in \cite{GB}.
 \bedefi \label{def:weakA-frame}
 Let $A$ be a  densely defined operator on $\H$.
 A \emph{(continuous) weak $A$-frame} is a function $\phi : x\in X \mapsto \phi_x$ such that, for all $u\in \D(A^*)$, the map  $x\mapsto \ip{u}{\phi_x}$ is a measurable function on $X$ and, for some $\alpha >0,$
\be \label{eq-weakAframe}
 \alpha \norm{}{A^\ast u }^2 \leq { \int_X |\ip{u}{\phi_x}|^2 \ud\mu(x) < \infty}, \qquad \forall \, u \in \D(A^\ast).
\en
   \findefi

 If $X= \NN$ and $\mu$ is the counting measure, we recover the discrete situation (so that the word `continuous' is superfluous in the definition above).
We get a  { simpler} situation   when $A$ is bounded and $\phi$  is Bessel.   { This is in fact  the construction of  G\u{a}vru\c{t}a \cite {gavruta}.}

Now we  introduce a structure that generalizes both concepts of   lower semi-frame and weak $A$-frame.
 We follow mostly the terminology of \cite{GB} and keep the term ``weak'' because the notion leads to a weak decomposition of the range of the operator $A$ (see Theorem \ref{th_char_continuous weak_A_frame}).

We begin with giving the following definitions.

\bedefi \label{defn_phiext}
Let $A$ be a  densely defined operator on $\H$,
		$\phi : x\in X \mapsto \phi_x$ a function  such that, for all $u\in \D(A^*)$, the map  $x\mapsto \ip{u}{\phi_x}$ is  measurable on $X$.
We say that a closed operator $B$ is a {\em $\phi$-extension} of $A$ if
$$ 
\mbox{$A\subset B$ \;and\; $\D(B^*)\subset \D(C_\phi)$.}
$$
We denote by ${\mc E}_\phi(A)$ the set of $\phi$-extensions of $A$.
\findefi
	{\berem It worths noting that if $A$ has a  $\phi$-extension, then $A$ is automatically closable.\enrem}
\bedefi
Let $A$ and $\phi$ be as in Definition \ref{defn_phiext}.
		Then $\phi$ is called a \emph{weak lower $A$-semi-frame}
		if $A$ admits a $\phi$-extension $B$ such that $\phi$ is a  weak $B$-frame.
		\findefi
Let us put $\D(A,\phi):=\D(A^*)\cap \D(C_\phi)$. If $\phi$ and $A$ are  as in Definition \ref{defn_phiext},
and $\displaystyle B:=\left(A^* \upharpoonright \D(A,\phi)\right)^*$ is
a $\phi$-extension of $A$,   it would be the smallest possible extension for which $\phi$ is a weak $B$-frame, but in general we could have a larger extension enjoying the same property.
Indeed, if $B$ is a closed extension of $A$ such that $\phi$ is a weak $B$-frame, we have
$$
A\subset A^{**} \subset \left( A^* \upharpoonright \D(A,\phi)\right)^* \subset B.
$$

\berems
\begin{enumerate}
\item
If $A$ is bounded, $\D(A^*)=\H$ and we recover the notion of lower semi-frame, under some minor restrictions on $A$,
hence the name (see Proposition
 \ref{prop: A-contr lsf is a lsf}).
	
\item If  $A$ is  a  densely defined operator on $\H$ such that the integral on the right hand side of \eqref {eq-weakAframe}
 is finite for every $f\in \D(A^*)$, then  {$\D(A^*)\subset\D(C_\phi)$ and} the
 weak lower  $A$-semi-frame $\phi$ is, in fact, a  weak $A$-frame, in the sense of  Definition \ref{def:weakA-frame}.

\item
Let us assume that $\phi$ is both a lower semi-frame  and a  weak  $A$-frame, then we have simultaneously 
\begin{align*}
{\sf m}  \norm{}{f}^2 &\leq  \int_{X}  |\ip{f}{\phi_{x}}| ^2 \, \ud \mu(x), \quad \forall \, f \in \H,
\\
\alpha \norm{}{A^\ast f }^2 &\leq \int_X |\ip{f}{\phi_x}|^2 \ud\mu(x) < \infty, \quad \forall \, f \in \D(A,\phi) =\D(A^\ast)\cap \D(C_\phi).
\end{align*}
It follows  that
\begin{equation}\label{eq-lax}\alpha'(\norm{}{f}^2  + \norm{}{A^\ast f }^2 ) \leq \int_X |\ip{f}{\phi_x}|^2 \ud\mu(x) < \infty, \quad \forall \,f \in \D(A^\ast)\cap\D(C_\phi)
\end{equation} {with $\alpha'\leq \frac12\min\{{\sf m},\alpha\}$.}
If we consider the domain $\D(A^*)$ with its graph norm
($  \|f\|_{A^\ast} = (\norm{}{f}^2  + \norm{}{A^\ast f }^2 )^{1/2}$, $f\in \D(A^\ast)$), we are led to the triplet of Hilbert spaces
$$
\H(A^\ast) \subset \H \subset \H(A^\ast)^\times,
$$
as discussed in Section \ref{sec2}.
Let us consider the  sesquilinear form $\Omega_\phi$ defined in \eqref{defn} and suppose in particular that $\D(A^\ast)=\D(\Omega_\phi)=\D(C_\phi)$. Then using Proposition \ref{prop: closedness C}, it is not difficult to prove that $\Omega_\phi$ is closed in $\H(A^\ast)$ and then bounded. Thus, there exists $\gamma >0$ such that, for every $f \in \D(A^\ast)$,
\begin{equation*}
\alpha'(\norm{}{f}^2  + \norm{}{A^\ast f }^2 ) \leq \int_X |\ip{f}{\phi_x}|^2 \ud\mu(x) \leq \gamma (\norm{}{f}^2  + \norm{}{A^\ast f }^2 ) .
\end{equation*}
One could notice that \eqref{eq-lax} is similar to  a frame condition.

The inequality \eqref{eq-lax} 
says that the sesquilinear form $\Omega_\phi$ defined in \eqref{defn} is coercive  
on $\H(A^\ast)$  and thus  the Lax-Milgram theorem applies \cite[VI \S 2, 2]{Kato} or \cite[Lemma 11.2]{Schmudg}. This means that for every $F \in \H(A^\ast)^\times$ there exists $w\in \H(A^\ast)$ such that
$$
\ip{F}{f} = \Omega_\phi (w, f) =\int_X\ip{w}{\phi_x}\ip{\phi_x}{f}\ud\mu(x), \quad \forall f \in \H(A^\ast).
$$
Therefore, in the case under consideration, we get expansions in terms of $\phi$ of elements that do not belong to the domain of $A^*$; in particular, those of $\H$. The price to pay is that the form of this expansion is necessarily {weak} since vectors of $\H$ do not belong to the domain of the analysis operator $C_\phi$.
\end{enumerate}
\enrems

  In the sequel we will need the following  

\begin{lemma}\cite[Lemma 3.8]{BC}
\label{doug gen unb}
	Let $(\H,\|\cdot \|),(\H_1,\|\cdot \|_1)$ and $(\H_2,\|\cdot \|_2)$ be Hilbert spaces and $T_1:\D(T_1)\subseteq \H_1\to \H$, $T_2:\D(T_2)\subseteq \H\to \H_2$ densely defined operators.
	Assume that
	$T_1$ is closed and
	$\D(T_1^*)=\D(T_2)$.\\ If
		 $\|T_1^* f\|_1\leq \lambda \|T_2f\|_2$ for all $f\in \D(T_1^*)$ and some $\lambda>0$, then
there exists a bounded operator $U\in \B(\H_1,\H_2)$ such that $T_1=T_2^* U$.
		\end{lemma}
	{
\berem\label{rem: doug gen unb}  Lemma \ref{doug gen unb} is still valid if  we replace closedness of $T_1$ by its closability, and in this hypothesis,  $\overline{T_1}=T_2^*U$.\enrem}

{In literature \cite{Z}, two measurable functions $\psi$ and $\phi$}  are said to be
\emph{dual} to each other if one has
\be\label{eq:dual}
\ip{f}{g}=\int_X\ip{f}{\phi_x}\ip{\psi_x}{g}\ud\mu(x),\qquad\forall f,g\in\H.
\en
If $\phi$ is a lower semi-frame of $\H$, then its dual $\psi$ is a Bessel mapping of $\H$ \cite{ACT}. In addition, if $\D(C_\phi)$ is dense, its dual $\psi$ is an upper semi-frame.

However this definition is too general, in the sense that the right hand side may diverge for arbitrary $f,g\in \H$. A more useful definition will be given below, namely \eqref{eq:dual1}.

A notion of duality related to a given operator $G$ can be formulated as follows.

	\bedefi\label{def24''}
 Let $G$ be a densely defined operator and 	$\phi : x\in X \mapsto \phi_x$ a function  such that, for all $u\in\D(G^*)$  the map  $x\mapsto \ip{u}{\phi_x}$ is a measurable function on $X$.
	Then a  function $\psi:x\in X\mapsto\psi_x\in\H$ 	  such that, for all $f\in\D(G)$  the map  $x\mapsto \ip{f}{\psi_x}$ is a measurable function on $X$ is called a {\it  weak $G$-dual of} $\phi$  if

	\be\label{eq24''}
	\ip{G f}{u}=\int_X  \langle f | \psi_x\rangle \ip{\phi_x }{u}\ud\mu(x),\quad \forall f\in\D(G)\cap\D(C_\psi), \forall u\in\D(G^*)\cap\D(C_\phi).
	\en
	\findefi 
This is a generalization of the notion of weak $G$-dual in \cite{GB}.

\berems \label{rem: on weak A dual 2}
\begin{itemize}
	\item[$(i)$] The  weak $G$-dual $\psi$ of  $\phi$ is not unique, in general. On the other hand, Definition \ref{def24''} could be meaningless. For instance, if either  $\D(G)\cap\D(C_\psi) =\{0\}$ or $\D(G^*)\cap\D(C_\phi)=\{0\}$, then everything is ``dual''.
	
	\item[$(ii)$] Note that, if $\phi$ is a weak $G$-frame, then there exists a weak $G$-dual $\psi$ of  $\phi$ such that relation \eqref {eq24''} must hold only for  $\forall f\in\D(G), \forall u\in\D(G^*)$ indeed $\D(G^*)\subset\D(C_\phi)$ and by Theorem 3.20 in \cite{GB} there exists a Bessel weak $G$-dual $\psi$ of $\phi$, hence $\D(G)\subset\D(C_\psi)=\H$.
\end{itemize}
\enrems

\beex
Given a densely defined operator $G$   on a separable Hilbert space $\H$, we   show two examples of $G$-duality (see  \cite[Ex. 3.10]{GB}).
\begin{enumerate}
\item[(i)] Let $(X,\mu)$ be a locally compact, $\sigma$-compact measure space and let $\{X_n\}_{n\in\mathbb{N}}$ be a covering of $X$ made up of countably many measurable disjoint sets of finite measure. Without loss of generality we suppose that $\mu(X_n)>0$ for every $n\in\mathbb{N}$. Let $\{e_n\}\subset\D(G)$ be an orthonormal basis of $\H$
and consider $\phi$, with $\phi_x= \frac{Ge_n}{\sqrt{\mu(X_n)}}$, $x\in X_n, \forall n\in\mathbb{N}$,  then $\phi$ is a  weak $G$-frame, see \cite[Example 3.10]{GB}. One can take $\psi$ with $\psi_x=\frac{e_n}{\sqrt{\mu(X_n)}}$, $x\in X_n, \forall n\in\mathbb{N}$.
\item[(ii)] If  $\phi:=G\zeta$, where $\zeta:x\in X\mapsto\zeta_x\in\D(G)\subset\H$ is a continuous frame for $\H$, then one can take as $\psi$ any  dual frame of $\zeta$.
\end{enumerate}
\enex

\section{Lower atomic systems}\label{subsec41}

 \begin{thm}\label{thm: A=RM supset DM} Let $(X,\mu)$ be a locally compact, $\sigma$-compact measure space, $A$ a densely defined operator and $\phi:x\in X\mapsto \phi_x\in\H$ a map such that, for every $u\in\D(A^*)$, the function $x \mapsto \ip{u}{\phi_x}$ is measurable on $X$.
	Then the following statements are equivalent.
	\begin{enumerate}
		\item[(i)] $\phi$ is  a weak lower $A$-semi-frame for $\H$.
		\item[(ii)] ${\mc E}_\phi(A) \neq \emptyset$ and for every $B \in {\mc E}_\phi(A)$,  there exists a closed densely defined  extension $R$ of $C_\phi^*$,   with  $\D(R^*)=\D(B^*)$, such that 
		$B$ can be decomposed as $B=RM$ for some $M\in \B(\H,L^2(X,\mu))$.
	\end{enumerate}
\end{thm}
\bdim 
 We proceed as in  \cite[Theor. 3.16]{GB}.
 \\
	\un{(i)$\Rightarrow$(ii)}	
If $\phi$ is  a weak lower $A$-semi-frame for $\H$, by definition, there exists $B\in {\mc E}_\phi(A)$.	{Consider  $E:\D(B^*)\to L^2(X,\mu)$ given by $(E u)(x) = \ip{u}{\phi_x }$, $\forall u\in\D(B^*)$, }$x\in X$  which is a restriction of the analysis operator $C_\phi $. $E$ is closable and densely defined.
	
	Apply Lemma \ref{doug gen unb}  to $T_1:=B$, 
 	 and $T_2:=E$,  noting that
	$\|Eu\|_2^2=\int_X	|\ip{u}{\phi_x }|^2\ud\mu(x)$, $u\in \D( B^*)$. Thus there exists $M\in \B(\H,L^2(X,\mu))$ such that    $B=E^*M$.  Then the statement is proved by taking $R=E^*$, indeed $R=E^* \supseteq C_\phi^*$ and  $\D(R)\supset\D(C_\phi^*)$ is dense because $C_\phi $ is closed and densely defined. Note that we have
	$\D(B^*)=\D(R^*)$; indeed $\D(R^*)=\D(\overline{E})$,
	$$
	\D(B^*)\subset\D(\overline{E})=\D(M^*\overline{E})\subset\D((E^*M)^*)=\D(B^*),
	$$
	hence in particular $E$ is closed, recalling that $\D(E)=\D(B^*)$.
	\\
	\un{(ii)$\Rightarrow$(i)}  Let $B\in {\mc E}_\phi(A)$; For every $u\in\D(B^*)=\D(R^*)$
	$$
	\|B^* u\|^2=\|M^*R^*u\|^2\leq\|M^*\|^2\|R^*u\|^2=\|M^*\|^2\int_X|\ip{u}{\phi_x}|^2\ud\mu(x)<\infty
	$$	since $R^*\subset C_\phi $. This proves that $\phi$ is a weak lower $A$-semi-frame.
	\edim

Generalizing the notion of continuous weak atomic system for $A$ \cite{GB}, we consider the following

\bedefi\label{def: continuous weak lower atomic system for A}
Let $A$ be a  densely defined operator on $\H$. A  \emph{  lower atomic system for $A$} is a function  $\phi:x\in X\mapsto\phi_x\in\H$ such that
 \begin{itemize}
\item[(i)] for all $u\in\D(A^*)$, the map $x \mapsto \ip{u}{\phi_x}$ is a measurable function on $X$;
\item[(ii)] the operator $A$ has a closed extension $B$ such that $\D(B^*)\subset \D(C_\phi)$; i.e., ${\mc E}_\phi(A)\neq \emptyset$;
\item[(iii)] there exists $\gamma>0$ such that, for every $f\in\D(A)$, there exists $a_f\in L^2(X,\mu)$, with
	$\|a_f\|_2=\left( \int_X |a_f(x)|^2\ud\mu(x)\right) ^{1/2}\leq \gamma\|f\|$ and
	\begin{equation*}
	{\ip{Af}{u}=\int_X  a_f(x) \ip{\phi_x }{u}\ud\mu(x), \qquad \forall u\in \D(B^*)}.
	\end{equation*}
\end{itemize}

\findefi
 We have  choosen not to call $\phi$ a  \emph{weak} lower atomic system for $A$ for brevity, even if it leads to a weak decomposition of the range of the operator $A$.
 \\

Theorem 3.20 of \cite{GB} gives a characterization of  weak atomic systems for $A$ and  weak $A$-frames. The next theorem yields the corresponding result for weak lower $A$-semi-frames.
\begin{thm}
	\label{th_char_continuous weak_A_frame}
Let $(X,\mu)$ be a locally compact, $\sigma$-compact measure space,   $A$  a  densely defined operator in  $\H$ and  $\phi:x\in X\mapsto \phi_x\in\H$ a function such that, for all $u\in\D(A^*)$, the map $x \mapsto \ip{u}{\phi_x}$ is  measurable on $X$. Then the following statements are equivalent.
	\begin{itemize}
		\item[(i)] $\phi$ is a lower    atomic system for $A$;
		\item[(ii)] $\phi$ is a weak lower $A$-semi-frame {for $\H$};
		\item[(iii)] ${\mc E}_\phi(A) \neq \emptyset$ and, for every $B\in {\mc E}_\phi(A)$,
		  $\phi$ has a Bessel   weak {$B$-dual} $\psi$.
	\end{itemize}
\end{thm}

\bdim \un{(i)$\Rightarrow$(ii)}
\\	{Consider a $\phi$-extension $B$ of $A$.}	By the density of $\D(A)$,   we have, for every    $u\in\D(B^*)$ 
	\begin{align*}\|B^* u\|&=\sup_{f\in \H,\|f\|=1}\left|\ip{B^*
			u}{f}\right|	=\sup_{f\in \D(A),\|f\|=1}\left|\ip{B^*
			u}{f}\right|\\	&=
		\sup_{f\in\D(A),\|f\|=1}|\ip{u}{Bf}| =	\sup_{f\in\D(A),\|f\|=1}|\ip{u}{Af}|
		\\&=\sup_{f\in \D(A),\|f\|=1}\left|\int_X
		\overline{a_f(x)}\ip{u}{ \phi_x }\ud\mu(x)\right|\\ &\leq
		\sup_{f\in \D(A),\|f\|=1}\left(\int_X
		|a_f(x)|^2\ud\mu(x)\right)^{1/2}\left(\int_X |\ip{u}{
			\phi_x }|^2\ud\mu(x)\right)^{1/2}\\&\leq \gamma\left(\int_X |\ip{u}{
			\phi_x }|^2\ud\mu(x)\right)^{1/2}<\infty\end{align*}
	for some $\gamma>0$; 
	 the last but one inequality is due to the fact that $\phi$ is a lower atomic system for $A$ and the last one to the inclusion $\D(B^*)\subset\D(C_\phi)$. Then, $\phi$ is a weak lower $A$-semi-frame. \\

	\un{(ii)$\Rightarrow$(iii)}
\\Following the proof of Theorem \ref{thm: A=RM supset DM}, 
for every  $\phi$-extension  $B$ of $A$ there exists a closed densely defined extension $R$ of $C_\phi^*$, with $\D(R^*)=\D(B^*)$,
	such that $B=RM$ for some $M\in \B(\H,L^2(X,\mu))$.

	By the Riesz representation theorem, for every $x\in X$ there exists a unique vector $\psi_x\in\H$ such that $(Mh)(x)=\ip{h}{\psi_x}$, for every $h\in \H$. The function $\psi:x\in X\mapsto\psi_x\in\H$ is Bessel. 
	 Indeed,  	
	\begin{eqnarray*} 
	\int_X
		|\ip{h}{\psi_x}|^2 \ud\mu(x)& {=}&  \int_X |(Mh)(x)|^2\ud\mu(x)\\&=&\nonumber\|Mh\|_2^2\leq\|M\|^2\|h\|^2, \qquad\forall h\in\H.
	\end{eqnarray*}	Hence $\D(C_\psi)=\H$.
	Moreover,  for $f\in \D(B){\cap\D(C_\psi)=\D(B)}$, $u\in \D(B^*)=\D(R^*)\subset\D(C_\phi)$
	\begin{eqnarray*}
		\ip{Bf}{u}&=&	\ip{R M f}{u}=\ip{M f}{R^* u}_2\\
		&=&\int_X\ip{f}{\psi_x}\ip{\phi_x}{u}\ud\mu(x).
	\end{eqnarray*}

	\un{(iii)$\Rightarrow$(i)}
	\\ It suffices to take, {for every fixed $\phi$-extension  $B$ of $A$}, $a_f:x\in X\mapsto a_x(f)=\ip{f}{\psi_x}\in\mathbb{C}$ for all $f\in \D(B)$. Indeed, $a_f\in L^2(X,\mu)$ and, for some $\gamma> 0$, we have $\int_{X} |a_x(f)|^2 \ud\mu(x) = \int_{X}|\ip{f}{\psi_x} |^2 \ud\mu(x) \leq \gamma  \|f\|^2$, since $\psi$ is a Bessel function. Moreover, by definition of weak $B$-dual, we have  $\ip{Bf}{u}=\int_X  a_f(x) \ip{\phi_x }{u}\ud\mu(x)$, for $f\in \D(C_\psi)\bigcap\D(B)=\D(B), u\in \D(B^*)\subset \D(C_\phi)$. Indeed we note that $\D(C_\psi)=\H$ since $\psi$ is a Bessel function.
\edim

\berem 
We don't know if $ \psi$  is a weak upper $A$-semi-frame, in the sense of  Definition
	\ref{def:weak A-upper semi-frame}, indeed $\psi$  needs not to be $\mu$-total, that is, $\int_X	|\ip{f}{\psi_x}|^2\neq 0$ for every $f\in\H$, $f\neq 0$. 
	\enrem

\section{Duality and Weak upper $A$-semi-frames} \label{subsec-duality}

If $C \in GL(\H)$, a \emph{frame controlled by the operator}  $C$ or $C$- \emph{controlled frame} \cite{contr} is a family of vectors $\phi = \left( \phi_n \in \Hil : n \in \Gamma \right)$,
	such that there exist two constants ${\sf m}_A>0 $ and ${\sf M}_A <\infty$ satisfying
{
\begin{equation}\label{eq:contrframineq1}
{\sf m_A} \norm{}{f}^2 \leq \sum \limits_n \ip{f }{ \phi_n}\ip{ C \phi_n}{ f }\leq {\sf M_A}  \norm{}{f}^2 ,
	{\forall}\,  f \in \H
\end{equation}
 or, to put it in a continuous form:
 \be	 \label{eq:contrframineq2}
	 {\sf m_A}  \norm{}{f}^2  \leq   \int_{X}  \ip{f}{\phi_x}\, \ip{C\phi_x}{f} \ud \mu(x)
	 \leq {\sf M_A}  \norm{}{f}^2 ,  \;  {\forall}\,  f \in \H.
	\en
	}

 According to Proposition 3.2 of  \cite{contr}, an    $A$-controlled frame is in fact a classical frame when the controlling operator belongs to $GL(\H)$. A similar result holds true for a
 	 weak lower $A$-semi-frame if
 $A$ is bounded as we show in Proposition
  \ref{prop: A-contr lsf is a lsf}. From there it follows that, if $A$ is bounded,  a weak lower $A$-semi-frame  has an upper semi-frame dual to it.

\berem\label{rem: rudin} We recall that a bounded operator $A$ is surjective if and only if $A^*$ is injective and  $\mathcal{R}(A^*)$ is norm closed (if and only if $A^*$ is injective and  $\mathcal{R}(A)$ is closed)
\cite[Theor. 4.14 and 4.15]{Rudin FA}.
\enrem

\begin{prop}\label{prop: A-contr lsf is a lsf}
 {Let $ A\in\B(\H)$ and  $\phi$ be a
	weak lower $A$-semi-frame. Assume that anyone of the following assumptions is satisfied:
	\begin{itemize}
		\item[(i)] $A^*$ injective, with $\mathcal{R}(A^*)$ 
		norm closed or
		\item[(ii)] $A^*$ injective, with $\mathcal{R}(A)$ closed or \item[(iii)]$A$ surjective.
	\end{itemize} }
	Then 	\begin{itemize}\item[(a)]$\phi$ is a 	lower semi-frame    of $\H$ in the sense of \eqref{eq:lowersf},\item[(b)] there exists an upper semi-frame $\psi$ dual to $\phi$.\end{itemize}
\end{prop}
\bdim
	(a) By {Remark \ref{rem: rudin} it suffices to prove (iii)}. By Theorem 4.15 in \cite{Rudin FA}, $A$ is surjective if and only if there exists $\gamma>0$ such that $\|A^*f\|\geq\gamma \|f\|$, for every $f\in \H$, then
	$$
	\gamma^2\alpha\|f\|^2\leq\alpha \|A^* f\|^2\leq\int_X|\ip{f}{\phi_x }|^2\ud\mu(x),\qquad \forall \,f\in\H.
	$$
	(b) The thesis follows from (a) and Proposition 2.1 (ii) in \cite{AT} (with $\{e_n\}$ an ONB of $\H$).
\edim

As explained above, the notion of duality given in \eqref{eq:dual} is too general. Therefore,
	in what follows  $\psi$ will be said to be \emph{dual to $\phi$} if one has
	\be \label{eq:dual1}
	\ip{f}{g}=\int_X\ip{f}{\phi_x}\ip{\psi_x}{g}\ud \mu(x),\qquad \forall f\in \D(C_\phi),\, g\in \D(C_\psi).
	\en

An interesting question is to identify a weak $A$-dual of a weak lower $A$-semi-frame. We expect one should generalize to the present situation the notion of upper semi-frame.
 We first  consider the next definition and examine its consequences.
	\bedefi \label{def:weak A-upper semi-frame}
Let $A$ be a  densely defined operator on $\H$. A  \emph{weak upper  $A$-semi-frame}
	for $\H$ is a function $\psi:x\in X\mapsto\psi_x\in\H$ such that, for all $f\in \D(A)$, the map  $x\mapsto \ip{f}{\psi_x}$ is  measurable on $X$ and there exists a closed extension $F$ of $A$
	and a constant $\alpha>0$ such that
	\begin{equation}\label{eq-contreusf}
	\int_X
	|\ip{u}{\psi_x }|^2\ud\mu(x)\leq\alpha \|F^* u\|^2,\qquad \forall \,u\in\D(F^*).
	\end{equation}
\findefi

\berems\begin{itemize}
	\item[(i)]From Definition \ref{def:weak A-upper semi-frame}  it is clear that $\D(F^*)\subset\D(C_\psi)$. 
	\item[(ii)] If $A\in\B(\H)$, then $\psi$ it is clearly a Bessel family.
\end{itemize}
\enrems

{\begin{cor}\label{cor: existence of an upper semiframe }
		Let $\psi$ be Bessel mapping of $\H$, and $ A\in\B(\H)$.  Assume that anyone of the following statements is satisfied:
		\begin{itemize}
			\item[(i)] $A^*$ injective, with 
$\mathcal{R}(A^*)$
 norm closed or
			\item[(ii)] $A^*$ injective, with $\mathcal{R}(A)$ closed or \item[(iii)]$A$ surjective.
		\end{itemize}
		Then $\psi$ is a weak upper $A$-semi-frame.
	\end{cor}
	\bdim By Remark \ref{rem: rudin} it suffices to prove (iii). {By Theorem 4.15 in \cite{Rudin FA},} we have just to note that
		$$\int_X|\ip{f}{\psi_x }|^2\ud\mu(x)\leq\gamma\|f\|^2\leq\alpha^2\gamma \|A^* f\|^2,\qquad \forall \,f\in\H.$$
	\edim
	
	\berem
		The previous result is true \emph {a fortiori} if $\psi$ is an upper semi-frame of $\H$.
	\enrem

	Summarizing Proposition \ref{prop: A-contr lsf is a lsf}, Corollary \ref{cor: existence of an upper semiframe } together with the preceding results  we have that \begin{cor}
		Let $A\in\B(\H)$. Assume that anyone of the following assumptions is satisfied:
		\begin{itemize}
			\item[(i)] $A^*$ injective, with $\mathcal{R}(A^*)$
			norm closed or
			\item[(ii)] $A^*$ injective, with $\mathcal{R}(A)$ closed or \item[(iii)]$A$ surjective.
		\end{itemize} and {let  $\phi$ be} a  weak lower $A$-semi-frame.
		Then there exists a 	weak upper $A$-semi-frame $\psi$   dual to $\phi$.
	\end{cor}
}

\begin{thm}\label{thm: A=RM supset DM upper} Let $(X,\mu)$ be a locally compact, $\sigma$-compact measure space, $A$ a densely defined operator and $\psi:x\in X\mapsto \psi_x\in\H$ a map such that, for every $f\in\D(A)$, the function $x \mapsto \ip{f}{\psi_x}$ is  measurable on $X$.
	Then the following statements are equivalent.
	\begin{enumerate}
		\item[(i)] $\psi$ is  a weak upper $A$-semi-frame for $\H$.
		\item[(ii)] For every  closed, densely defined   extension  $F$ of $A$ such that \eqref{eq-contreusf} holds true, there exists  a closed, densely defined  extension $Q$
		of $C_\psi^*$ such that $Q=FN$ for some  $N\in\B(L^2(X,\mu),\H)$.
	\end{enumerate}
\end{thm}
\bdim
	\un{(i)$\Rightarrow$(ii)}
	Let $\psi$ be  a weak upper $A$-semi-frame, then for every  closed extension $F$ of $A$
	for which \eqref{eq-contreusf} holds true,   consider the operator $E=C_\psi\upharpoonright\D(F^*)$. It is   densely defined, closable since $C_\psi$ is closed.
 Define an operator $O$ on $R(F^*)\subseteq \H $ as $OF^*f=E f\in L^2(X,\mu)$. Then $O$ is a well-defined bounded operator by \eqref{eq-contreusf}. Now we extend $O$ to the closure of $R(F^*)$ by continuity and define it to be zero on $R(F^*)^\perp$. Therefore $O\in \B(\H,L^2(X,\mu))$ and $OF^*= E $, i.e. $E^*=F O^*$ and the statement is proved by taking $Q=E^*$ and $N=O^*$.	
	\\		\un{(ii)$\Rightarrow$(i)}
	From $Q=FN$, with $Q$ a densely defined closed extension of $C_\psi^*$, we have that $Q^*=N^*F^*\subset C_\psi$. For every $u\in\D(F^*)=\D(N^*F^*)=\D((FN)^*)\subset\D(C_\psi)$		
	$$\|C_\psi u\|_2^2=\int_X|\ip{u}{\psi_x}|^2\ud\mu(x)=\|N^*F^* u\|_2^2\leq\alpha\|F^* u\|^2$$ for some $\alpha>0$.
	\edim

We  can now prove the following duality result, which suggests that Definition \ref{def:weak A-upper semi-frame} is convenient in this context.

\beprop \label{prop59}
Let  $A$ be a densely defined operator and $\psi$  a weak  upper  $A$-semi-frame.  
Let $F$ be a closed extension of $A$ satisfying \eqref{eq-contreusf} for some $\alpha >0$.
Assume that $\phi  \subset\D(A)$ 
is a weak  $F$-dual  of $\psi$ such that
\begin{itemize}
	\item[(a)] $F^*\D(F^*)\subset \D(C_\phi)$;
	\item[(b)] the function $x\to \|A\phi_x\|$ is in $L^2(X, \mu).$ 
\end{itemize}
Then $F\in {\mc E}_{A\phi}(A)$ and $A\phi$  is a weak lower $A$-semi-frame  with $F$ as $(A\phi)$-extension and lower bound
$\alpha^{-1}$; i.e.,
\begin{equation}\label{eq: Aphi lower A semi-frames}
\alpha^{-1}\, \|F^* u\|^2\leq\int_X
|\ip{u}{ { A\phi_x} }|^2\ud\mu(x),\qquad \forall \,u\in{ \D(F^*)\cap\D(C_\phi)}.
\end{equation}
\enprop
\bdim
	For every { $u\in\D(FF^*)$} 
	\begin{eqnarray*}\|F^*  u\|^2
		&=&\ip{F^*u}{F^*u}=\ip{ {F}F^*u}{u}
		\\	
		&=&		\int_X \ip{F^*u}{\phi_x } \ip{ \psi_x}{u} \ud\mu(x), \;\; \mbox{by {weak} $F$-duality}
		\\
		&\leq&		\left(\int_X|\ip{u}{\psi_x}|^2\ud\mu(x)\right)^{1/2}\left(\int_X |\ip{ F^*u
		}{
			\phi_x }|^2\ud\mu(x)\right)^{1/2}\\
		&\leq& \alpha^{1/2} \, \|F^*u\|\left(\int_X |\ip{ F^*u}{
			\phi_x }|^2\ud\mu(x)\right)^{1/2}\end{eqnarray*}
	{The right hand side of the previous inequality is finite because of (a).}
	
	Hence, $$\|F^*  u\|
	\leq \alpha^{1/2} \left(\int_X |\ip{u}{{ A
			\phi_x }}|^2\ud\mu(x)\right)^{1/2}, \quad \forall u\in \D(FF^*).$$

	Now we take into account that $\D(FF^*)$ is a core for $F^*$ by von Neumann theorem \cite[Theorem 3.24]{Kato}. Therefore, for every $u\in \D(F^*)$, there exists a sequence $\{u_n\}\subset \D(FF^*)$ such that $\|u_n -u\|\to 0$ and $\|F^*u_n -F^*u\| \to 0$. This implies, of course, that $\ip{ F^*u_n} 
		{\phi_x } \to \ip{ F^*u} 
		{\phi_x }$, for every $x\in X$. 
	Moreover, since $\{u_n\}$ is bounded, we have
	$$
	|\ip{F^*u_n}{\phi_x}|=|\ip{u_n}{F\phi_x}|\leq M \|F\phi_x\|,
	$$ 
	  for some $M>0$ and for every $x\in X$. 
	The assumption that $x\to \|A\phi_x\|$ is in $L^2(X, \mu)$ allows us to apply the dominated convergence theorem and conclude that
	$$
	\|F^*  u\| \leq \alpha^{1/2} \left(\int_X |\ip{u}{{ A\phi_x }}|^2\ud\mu(x)\right)^{1/2}, \quad \forall u\in \D(F^*).
			$$
	The right hand side of the latter inequality is finite again by (a),  hence $\D(F^*)\subset\D(C_{A\phi})$. This fact also implies that $F\in {\mc E}_{A\phi}(A)$ since, if $u\in \D(F^*)$, we get 
	
	$$
	\int_X |\ip{u}{{ A
			\phi_x }}|^2\ud\mu(x) = \int_X |\ip{F^*u}{{ 
			\phi_x }}|^2\ud\mu(x)= \|C_{A\phi} u\|^2<\infty.
			$$
\edim

\berems
\bee

\item
Note \eqref{eq: Aphi lower A semi-frames} can obviously be also written
$$
\alpha^{-1}\, \|h\|^2\leq\int_X
|\ip{h}{ \phi_x }|^2\ud\mu(x),\qquad \forall \,h\in\mathcal{R}(F^*).
$$
\item

For every {$f\in\D(F^*)$, with our new definition} by 
\begin{equation*}
\alpha^{-1}\int_X
|\ip{f}{\psi_x }|^2\ud\mu(x)\leq\, \|A^* f\|^2
\leq \alpha\int_X |\ip{f}{ A
	\phi_x }|^2\ud\mu(x)
	\end{equation*}
	it follows that $\|C_\psi f\|\leq\alpha\|C_{A\phi} f\| $ for every {$f\in\D(F^*)$}. Since $C_\psi$ is closed then
	$\D(C_\psi^*)$ is dense and \eqref{eq-contreusf} and \eqref{eq: Aphi lower A semi-frames}
imply that $\D(C_{A\phi})\subseteq{\D(F^*)}\subseteq\D(C_\psi)$ hence the last is dense too. 

\ene
\enrems

{ Another possibility is to mimic the notion of controlled frame \eqref{eq:contrframineq1} or  \eqref{eq:contrframineq2}, introduced  in \cite[Definition 3.1]{contr}.
 Because the operator $A$ is supposed to belong to $GL(\Hil)$, we end up with a generalized frame (and actually a genuine frame). It would be interesting to extend the definition  to an unbounded operator or at least to operators  less regular than elements of $GL(\H)$. 

  A possibility  is  to investigate the following generalization. Let $B$ be a linear operator with domain $\D(B)$.   Suppose that $\psi_n\in  \D(B)$ for all $n$.
Put
\begin{equation*} 
 \Omega_B(f,g) =\sum_n \ip{f}{\psi_n}\ip{B\psi_n}{g}, \forall\;f,g \in  \D( \Omega_B),
\end{equation*}
 where $\D(\Omega_B)$ is  some  domain of the sesquilinear form defined formally on the rhs. Following
 \cite[Sec. 4]{corso}, we may consider the form $\Omega_B$ as the form generated by two sequences, $\{\psi_n\}$ and $\{B\psi_n\}$. Then the operator associated to the form $ \Omega_B$ is precisely $B$, since one has
 $\ip{B f}{g} =  \Omega_B(f,g)$.
\\

A continuous version of \eqref{eq:contrframineq1}  would be
$$
 {\sf m}_{A} \norm{}{f}^2 \leq \int_{X}  \ip{f}{\psi_x}   \ip{A  \psi_x} {f } \ud \mu(x)
\leq { \sf M}_{A} \norm{}{f}^2 ,
\quad \mbox{for all}\;  f \in \Hil.
$$
and the sesquilinear form becomes
$$
 \Omega_A(f,g) = \int_{X}  \ip{f}{\psi_x}   \ip{A \psi_x} {g } \ud \mu(x)
\leq { \sf M}_{A} \norm{}{f}^2 , \quad \mbox{for all}\;  f, g\in \D(\Omega_A).
$$

From the last relation, we might infer two {alternative} possible definitions of  an \emph{upper $A$-semi-frame}, namely:
\bea
    &&\nonumber\int_{X}  |\ip{Af}{\psi_{x}}| ^2 \, \ud \mu(x)   \leq { \sf M}  \norm{}{f}^2 , \quad \forall \, f \in \D(A),  
\\
  &&  \int_{X}  \ip{f}{\psi_x}\, \ip{ \psi_x}{Af} \ud \mu(x)  \leq {\sf M} \norm{}{f}^2 ,  \quad \forall \, f \in \D(A),\label{eq56}
\ena
Actually the   definition \eqref {eq56} leads to that of an \emph{$A$-Bessel map}, provided that $\psi_x \in \D(A^*)$,
for all $x\in X$:
$$
\int_{X}  \ip{f}{\psi_x}\, \ip{ A^*\psi_x}{f} \ud \mu(x)  \leq {\sf M} \norm{}{f}^2 ,  \quad \forall \, f \in \D(A).
$$}

Further study will hopefully reveal which of the three definitions of  an  {upper $A$-semi-frame} is the most natural one.

\section{Examples}\label{sec-examples}

\subsection*{\small(1) A reproducing kernel \hs}

We start from the example of a lower semi-frame in a  reproducing kernel \hs\ described in increasing generality in \cite{at-reprodpairs2, ACT}.
Let $\H_K$ be  a reproducing kernel Hilbert space   of (nice) functions on
a measure space $(X, \mu)$, with kernel function $k_x, x\in X$, that is, $f(x)=\ip{f}{k_x}_K,\, \forall f\in\H_{K}$.
Choose a (real valued, measurable)  weight  function $m(x) >1$  and consider the unbounded self-adjoint  multiplication operator $(Mf)(x) = m(x) f(x),    \forall x\in X$, with dense domain $\D(M)$.  
For each $ n \in\NN$, define $\H_{n}= \D(M^n)$,  equipped with its graph norm,  and $\H_{\ov n}:= \H_{-n} =\H_{n}^{\times}$ (conjugate dual). Then we have the Hilbert scale $\{\H_n, \, n\in \ZN\}$: 
$$
\ldots \H_n \subset  \ldots \subset   \H_2  \subset   \H_1  \subset  \H_0 = \H_K\ \subset
   \H_{\ov 1}   \subset \H_{ \ov 2} \ldots \ \subset \H_{ \ov n}  \ldots 
$$
As an operator on the scale, which is a partial inner product space \cite{pip-book},  the operator $M$ has continuous representatives  $M_{n+1} \to M_{n}, n\in \ZN.$

Fix some $n > 1$ and define
  the measurable functions $\phi_x = k_x m^n(x), \psi_x = k_x m^{-n}(x)$, for every $x\in X$.. Then $\psi_x \in \H_{n}$\,,  for every $x\in X$,    and $\psi$  is an upper semi-frame,   whereas $\phi_x \in \H_{\ov n}$\,, for every $x\in X$, and $\phi$ is a lower  semi-frame. Also
 $C_\psi : \H_K \to   \H_n,  \, C_\phi : \H_K \to  \H_{\ov n}$  continuously.  
 One has indeed, for every $g\in\H_K,\ip{\psi_x}{g}_K = \ov{g(x)}\,m^{-n}(x) \in \H_{n}$  and  
 $ \ip{\phi_x}{g}_K = \ov{g(x)}\,m^n(x) \in \H_{\ov n}$.

Next choose a real valued, measurable function $x\mapsto a(x)$ 
such that $a(x) \leq m^n(x), \forall x\in X$, and define $A= A\ha$ as the multiplication operator by $a : (Af)(x) = a(x)f(x),\forall \,x\in X$.   
 Let $\D(A) = \H_{n}$.  Then $A\in \B(\H_{n})$ since $\|af\|\leq\|m^nf\|<\infty$, for every $f\in \H_{n}$ and since $a(x)m^{-n}(x)<1$ for every $x\in X$ and for every $f\in\H_{n}$, then $\mathcal{R}(A)\subset\D(M^{n})
=\H_n$. 
As an operator on the scale, $A$ has continuous representatives $A_{p, n+p}: \H_{n+p} \to \H_{p}$.

Then we have, $\forall f\in\D(A)= \H_{n}\subset \D(C_\phi),$
$$
	 \|Af\|^2  = \int_X   |f(x)|^2 \,a(x)^{2} \ud\mu(x)\leq \int_X   |f(x)|^2 \,m^{2n}(x) \ud\mu(x)
  = \int_X |\ip{f}{\phi_x }_K|^2\ud\mu(x) < \infty.
  $$
that is, $\phi$ is   a  weak $A$-frame for $\H_K$ .

The same holds for every self-adjoint operator $A'$ which is  the multiplication operator  by the measurable function $x\mapsto a'(x)$ such that $a'(x) \leq m^n(x), \forall x\in X$, and $\D(A')=\H_{n}$.

Let now the closed operator $B$ be a  $\phi$-extension of $A$, that is,
$$ A\subset B  \mbox{ and }\;   \H_n =  \D(A) \subset \D(B^*)\subset \D(C_\phi)  $$
and 
$$
\|B^*f\|^2  \leq  \int_X |\ip{f}{\phi_x }_K|^2\ud\mu(x) < \infty, \; \forall\, f\in \D(B^*).
$$
Then  $\phi$ is   a  weak lower $A$-semi-frame for $\H_K$.

\subsection*{\small(2) A discrete example}

A more general situation may be derived from the discrete example of Section 5.2  of \cite{at-reprodpairs2}.
Take a weight  sequence  $m:=\{|m_n|\}_{n\in\NN},  m_n\neq 0,$   where  $m\in\ell^\infty$  has a subsequence converging to zero (or $m\in c_0$).
Then consider the space $\ell^2_m$ with norm $\norm{\ell^2_m}{\xi}:=\sum_{n\in\NN}|m_n \xi_n |^2$. Thus we have the following triplet,  
  $$
 \ell^2_{1/m} \subset \ell^2  \subset \ell^2_m.
 $$
 Next,    for each $n\in\NN$,  define $\psi_n = m_n e_n$, where $ e:=\{e_n\}_{n\in\NN}$ is   an orthonormal basis in $\ell^2$. Then $\psi $  is an upper semi-frame and
   $C_\psi : \H \to    \ell^2_{1/m}$, continuously.
On the other hand,       $\phi:=\{(1/\ov{m_n})e_n\}_{n\in\NN}\}$ is a lower semi-frame and
$ C_\phi : \H \to  \ell^2_m$,   continuously.

In other words, $\psi = M e$ and $\phi = M^{-1} e$, where $M$ is the diagonal operator $M_n = m_n,   n\in\NN$.
In order to define a weak  lower $A$-semi-frame for $\ell^2$, we take another diagonal operator $A = \{a_n\} $ such that, for each $n\in\NN$ one has
$|a_n| \leq |m_n| ^{-1}$. Then, $\forall f\in\D(A),$
\beano
 \|Af\|^2  &=& \sum_{n\in\NN} |a_n|^2  |f_n|^2  = \sum_{n\in\NN} |a_n|^2 |\ip{f}{e_n}|^2 \leq \sum_{n\in\NN}|m_n| ^{-2} |\ip{f}{e_n}|^2
 \\
 &=& \sum_{n\in\NN}  |\ip{f}{\phi_n}|^2.
\enano
Thus $\phi$ is   a weak  $A$-frame for $\ell^2$.  
As in Example (1), we get a weak lower $A$-semi-frame for $\ell^2$ if we have a $\phi$-extension $B$ of $A$.

The same result holds true if one replaces the ONB $\{e_n\}$ by a frame $ \{\theta_n\}_{n\in\NN}$ :
$$
\alpha \norm{}{f}^2 \leq \sum_n |\ip{f}{\theta_n}|^2 \leq \beta \norm{}{f}^2, \, \forall f\in\H
$$ for some $\alpha,\beta>0$.
Since  { $m\in \ell^\infty$}, we can as well assume that $|m_n |< \delta, \, \forall n\in\NN$ for some $\delta>0$. Thus
$|1/\ov{m_n}|  >  1/\delta, \, \forall n$. Then for every $f,g\in\H$ we have
\beano
\sum_n |\ip{f}{\psi_n}|^2 &=\; { \sum_n |m_n|^2|\ip{f}{\theta_n}|^2 \leq \delta^2 \sum_n|\ip{f}{\theta_n}|^2 \leq \delta^2}\beta \,\norm{}{f}^2,
\\
\sum_n |\ip{g}{\phi_n}|^2 &=\;  
{ \sum_n \left|\dis\frac{1}{m_n}\right|^2|\ip{g}{\theta_n }|^2\geq \dis\frac{1}{\delta^2} \dis\sum_n \left| \ip{g}{\theta_n}\right|^2 \geq  \dis\frac{1}{\delta^2}}\; \alpha \norm{}{g}^2.
\enano
Thus, indeed, $\psi$ is an  upper semi-frame and $\phi$ is a  lower semi-frame. The rest of the construction follows.

\subsection*{ \small (3) A standard construction}

As explained in Section  \ref{sec2}, a standard construction of lower semi-frames stems from the consideration of a metric operator induced by an unbounded operator.

 Given a closed, densely defined, unbounded   operator $S$ with dense domain  $\D(S)$,  define the metric  operator $G= I+ S^\ast S$, which is unbounded  with bounded inverse.

Then, if we take  an ONB $\{e_n\}$ of  $ \D(G^{1/2}) =\D(S)$, contained in
$\D(S^\ast S) $, then $   \{\phi_n\}=\{G e_n\} =  \{( I+ S^\ast S)e_n\}$ is a lower semi-frame of $\H$  on $\D(S)$.

Now, if $A$ is a densely defined operator that satisfies the equation
$$
\alpha \norm{}{A^*f} \leq \norm{C_\phi}{f}, \;\forall \,f\in\D(A^*)
$$
instead of \eqref{eq:triple1}, then $\phi$ is a weak $A$-frame for $\H$.
As for the equivalent of \eqref{eq:triple2}, it is of course trivial.

 \subsection*{Acknowledgements} This work has been supported by the
Gruppo Nazionale per l'Analisi Matematica, la Probabilit\`{a} e le
loro Applicazioni (GNAMPA) of the Istituto Nazionale di Alta
Matematica (INdAM).

\vspace*{0.5cm}

\bibliographystyle{amsplain}

\end{document}